\documentclass[a4paper]{article}
\pdfoutput=1                         
\usepackage[english]{babel} 

\usepackage{csquotes}                
\usepackage{authblk}                 

\usepackage[
  backend=bibtex,
  sorting=none,
  style = numeric-comp,
]{biblatex}
\usepackage{amsmath}
\usepackage{amsfonts}
\usepackage{amssymb}

\usepackage{graphicx}                
\usepackage{wrapfig}

\usepackage{verbatim}
\usepackage{xcolor}
\usepackage{listings}
\usepackage{algorithm}    
\usepackage{algpseudocode}

\usepackage{caption}
\usepackage{tabularray}
\UseTblrLibrary{booktabs}
\usepackage{pgfplots}
\usepackage[colorlinks=true,linkcolor=blue,citecolor=green,urlcolor=blue]{hyperref}

\pgfplotsset{compat=1.17}            

\begin{document}

\selectlanguage{english}

\title{Optimization of Discrete Parameters Using the Adaptive Gradient Method and Directed Evolution}

\author{%
Andrei Beinarovich,
~~Sergey Stepanov\thanks{Corresponding author: \textless{}\href{mailto:steps@qudata.com}{steps@qudata.com}\textgreater{}},
~~Alexander Zaslavsky \\
\em{QuData AI, Dnipro, Ukraine}\\
\normalsize{\url{https://qudata.com/}} \\
}

\date{December 21, 2023} 

\large

\maketitle

\begin{abstract}
The problem is considered of optimizing discrete parameters in the presence of constraints.
We use the stochastic sigmoid with temperature and put forward the new adaptive gradient method CONGA.
The search for an optimal solution is carried out by a population of individuals.
Each of them varies according to gradients of the 'environment' and is characterized by two temperature parameters with different annealing schedules.
Unadapted individuals die, and optimal ones interbreed, the result is directed evolutionary dynamics.
The proposed method is illustrated using the well-known combinatorial problem for optimal packing of a backpack (0--1 KP).
\end{abstract}

\section{Introduction}

Deep learning models usually utilize real-valued parameters.
This is due to the fact that the only known optimization method applicable to models with millions or more parameters is gradient descent.
However, many tasks that machine learning is aiming to solve are intrinsically discrete in nature.
This discreteness may follow directly from the structure of the input or output data, and also can be optimal configuration for 
the latent feature space.

For example, natural language operates with tokens, the relationship between which can be described using tree structures.
A similar situation occurs when working with particles, molecules, graphs and other objects that are of a discrete nature.
Over the years, researchers have become adept at 'encoding' discrete objects using continuous parameters (using embedding vectors, 
for example). However, this accomplishment comes at the cost of introducing a large number of redundant parameters. 

There are numerous real world problems that require a solution of optimization task with
constrained parameters, some of the most prominent examples being the traveling salesman problem (TSP), 
ranking and sorting algorithms, routing optimization, the 0--1 knapsack problem, and many other. Application to these combinatorial 
problems of the traditional ML methods that rely on the gradient descent is not straightforward, as the solution space is discrete 
and objective functions are not differentiable, as they involve step-wise decisions or discrete choices, rendering gradients 
undefined or impractical to compute. Also traditional gradient-based methods may not naturally integrate constraint satisfaction 
into the optimization process. Researchers put forward a number of approaches to tackle these difficulties. After the papers \cite{chung2016hierarchical, maddison2017concrete, jang2017categorical} introduced the continuous relaxation of 
discrete random variables it became possible to use the gradient method for discrete parameters. 

In the present paper we put forward a new adaptive gradient method for discrete optimization tasks.
The layout of the paper is as follows. In Section~\ref{relatedwork} we give a brief overview of the research conducted
recently in the domain of discrete optimization problems with constraints. Section~\ref{method} presents our new method for
solving such problems, CONstrained Gradient descent with Adaptation (CONGA). The application of hot and stochastic sigmoids for 
gradient computation is discussed in detail; the optimization procedure is carried out by a population of individuals in directed 
dynamic evolution. Section~\ref{experiments} demonstrates the performance of the CONGA method applied to one of the benchmark
discrete problems, the 0--1 Knapsack Problem. The numerical details of the experiments and comparison of our results to 
branch and bound, genetic algorithm, greedy search algorithm and simulated annealing are compiled in the Appendices.

\section{Related Work}\label{relatedwork}

In \cite{chung2016hierarchical} a multi-scale approach was proposed, called the hierarchical multi-scale recurrent 
neural network, that can capture the latent hierarchical structure in the sequence by encoding the temporal dependencies
with different timescales. It found that the straight-through estimator is one efficient way of training a model 
containing discrete variables is the straight-through estimator, which was first introduced in~\cite{Bengio2013} and works 
by decomposing the operation of a binary stochastic neuron into a stochastic binary part and a smooth differentiable part, 
approximating the expected effect of the pure stochastic binary neuron.

Authors of~\cite{maddison2017concrete} introduced the continuous relaxation of discrete random variables, using Gumbel-Softmax trick. 
They introduce the Concrete distribution with a closed form density parameterized
by a positive temperature parameter. The essence of the approach is that the zero temperature limit of a Concrete distribution 
corresponds to a discrete distribution of parameters. Thus, optimization of an objective over an architecture with discrete
stochastic nodes can be accomplished by gradient descent on the samples of the corresponding continuous relaxation.

The concept of employing a Gumbel-Softmax as a relaxation technique for discrete random variables was also explored 
by~\cite{jang2017categorical}. In their approach, the relaxed objective does not incorporate density; rather, all 
elements of the graph are evaluated using the relaxed stochastic state of the graph, which includes discrete log-probability 
calculations. Gumbel-Softmax estimator outperformed gradient estimators on both Bernoulli variables and categorical variables; 
the most important result is a differentiable approximate sampling mechanism for
categorical variables that can be integrated into neural networks and trained using standard backpropagation.
A recent survey of the applications of Gumbert-softmax trick to numerous machine learning tasks, as well as its extensions 
and modifications is available in~\cite{huijben2022review}.

There are a wide range of machine learning problems and applications, where discrete optimization may be profitably applied. 
We will outline here some interesting examples of recent research in the field.
Vector quantized variational autoencoders (VAE) were first introduced in~\cite{van2017neural}; the approach uses discrete rather 
than continuous latent space of categorical variables, which is arguably better corresponds to many language tasks. 
Another discretization of the VAE was considered in~\cite{ramesh2021zero} for the task of text-to-image generation 
(the first incarnation of DALL-E); this paper implements the autoregressive transformer model and utilizes 
the Gumbel-Softmax optimization.

In the paper~\cite{bao2021beit} authors proposed using VAE for a masked image modeling task to pretrain vision transformers (BEiT) 
in a self-supervised manner; this self-supervised model learned to segment semantic regions and object boundaries 
without prior human annotation.

Another productive family of models for the discrete optimization are graph neural networks. MolGAN~\cite{decao2022molgan} 
is a generative model for small molecular graphs that circumvents the need for expensive graph matching procedures or node 
ordering heuristics, adapting instead  generative adversarial networks (GANs) to operate directly on graph-structured
data, in combination with a reinforcement learning, and generate molecules with specific desired chemical properties.

The medical Generative Adversarial Network (medGAN) introduced in~\cite{Choi2017} can generate realistic synthetic high-dimensional patient records, containing  discrete variables (e.g., binary and count features) via a combination of an autoencoder and generative adversarial networks.

Sorting and ranking of the input objects is a key step in many algorithms and an important task in itself (e.g. ranking the documents 
in web or database search). 
The Gumbel-Sinkhorn algorithm, an analog of Gumbel-Softmax for permutations, was constructed in~\cite{mena2018learning} and applied to
several diverse problems (number sorting, jigsaw puzzle solving, identifying C. elegans neural signals), outperforming state-of-the-art 
neural network baselines.
Recently proposed Sparse Sinkhorn Attention~\cite{Tay2020} is a new memory-efficient method for learning to attend,
based on differentiable sorting of internal representations for transformer networks. This approach learns to generate 
latent permutations over sequences, and with sorted sequences, is able to compute global attention with only local windows,
mitigating the memory constraints.

Let's us now turn to 0--1 Knapsack Problem (KP)~\cite{Kellerer2004}, a classic NP-hard problem that is frequently encountered 
in real-world scenarios, including decision-making processes, resource allocation, cryptography, and computer vision. 
It involves selecting the best combination of items to maximize an objective function while satisfying a certain budget constraint. 
Over time, numerous algorithms have been devised to efficiently solve the KP. While some rely on traditional combinatorial 
optimization techniques, others harness the power of deep learning methods. Traditional methods such as integer linear programming 
and dynamic programming provide exact solutions but often suffer from scalability issues. As the dimensionality of the problem increases, 
the runtime of classical algorithms tends to grow exponentially, which can be a significant drawback in practical applications.
Even for the relatively low-dimensional KP there are classes of instances that pose serious difficulties for the exact methods,
as was shown in~\cite{pisinger2005hard}. Here, several instance dataset generation algorithms were constructed, which were later used 
by many researchers to produce the benchmarks for consistent comparison of the results. 
Several new computationally challenging problems were recently described in~\cite{Jooken2023}.

A comparative study of meta-heuristic optimization algorithms for 0--1 KP is given in~\cite{Ezugwu2019,Laabadi2018,Puchinger2010}, 
which provide a theoretical review and experimental results of applying a number of methods, including genetic algorithms, 
simulated annealing, branch and bound, dynamic programming, and greedy search algorithm. 
In recent years, a number of attempts were made at solving the combinatorial optimization problems by deep learning methods. 
For example, the traveling salesman problem is being solved by recurrent neural networks in~\cite{Vinyals2015}. 
A machine learning approach to getting an approximate solution to the KP is proposed in the paper~\cite{Gu2018}, 
consisting in randomly generating of large quantity of samples of 0--1 KP and applying supervised learning to train
a pointer network to this task. 
Another novel method to solve the large-scale hard KPs using advanced deep learning techniques was proposed in~\cite{li2021novel}.
The authors designed an adaptive gradient descent method that proved able to optimize the KP objective function, and
conducted extensive experiments with state-of-the-art baseline methods and benchmark KP datasets.

\newpage
\section{Proposed Method}\label{method}

Let's consider the problem of maximizing the function $v(\mathbf{x}) > 0$ of binary variables $\mathbf{x}=\{x_1,...,x_n\}$ 
in the presence of a constraint:
\begin{equation}
v(\mathbf{x}) = \max,~~~~~~w(\mathbf{x}) \le 0,~~~~~~x_k \in \{0,1\}.
\end{equation}
Using the gradient method $\mathbf{x}\mapsto \mathbf{x} - \lambda \nabla L$ with learning rate $\lambda$, 
we will look for the minimum of the following loss function:
\begin{equation}\label{loss}
L(\mathbf{x}) = -v(\mathbf{x}) + \frac{\gamma}{\nu}\,\bigr[\max(0, w(\mathbf{x}))\bigr]^\nu,
\end{equation}
where $\gamma \ge 0$, $\nu > 0$ are hyperparameters that determine the magnitude of the penalty for 
violating the constraint $w(\mathbf{x})\le 0$.
In the paper~\cite{li2021novel}, a 'game-theoretic approach' was used to justify the choice of $\nu=2$,
which resulted from the use of the quadratic regularization term. In general, the choice of $\nu$ is arbitrary and depends 
on the nature of the function $w(\mathbf{x})$.
For smooth functions, it is reasonable to set $\nu > 1$ so that there is no gradient jump at the transition point from 
the allowed region $w(\mathbf{x})\ge 0$ to the forbidden one $w(\mathbf{x}) < 0$ and back.

\subsection{The Hot Sigmoid}

We cannot directly minimize the loss function $L(\mathbf{x})$ using the gradient method because the derivatives by the integer-valued 
variables $x_k$ are not defined.

\begin{wrapfigure}{r}{0.5\textwidth}\centering
    \includegraphics[width=0.45\textwidth]{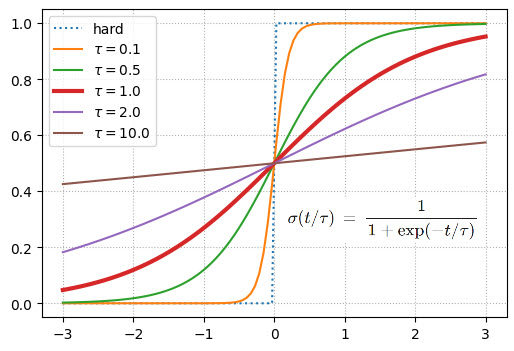}    
    \caption{\small Heaviside function (dashed) and sigmoid for different values of temperature $\tau$}\label{sigmoid}
\end{wrapfigure}   
We thus will compute the error gradient using the real-valued variable $t\in [-\infty, \infty]$ (the last letter in the word \emph{logit}).
The Heaviside function $x=H(t)$ of the logit is equal to $x=1$ (for $t > 0$) and
$x=0$ otherwise (in the graph \ref{sigmoid} $H(t)$ is shown by the dashed line, while the logit $t$ is plotted on the horizontal axis).

After substituting $x_k = H(t_k)$ into $L(\mathbf{x})$, we see the loss now depends on the continuous parameters $t = \{t_1,...,t_n\}$, 
still only taking on the 'allowed values' since $x_k$ are still discrete quantities.
However, the derivative of $H(t)$ at $t\neq 0$ is zero. Hence, the error gradient will also vanish.

\begin{wrapfigure}{r}{0.5\textwidth}\centering
    \includegraphics[width=0.45\textwidth]{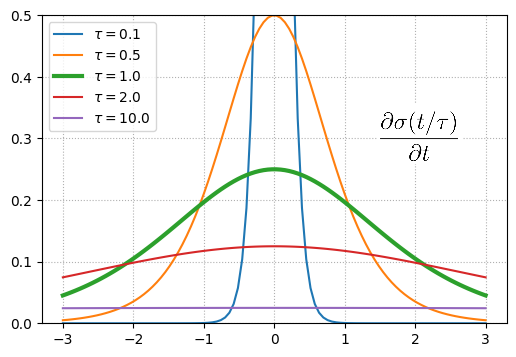}    
    \caption{\small The derivative of the sigmoid with temperature at different values of $\tau$}\label{sigmoid-df}
\end{wrapfigure}
To solve this problem,
following the \cite{chung2016hierarchical}, we will do the following. 
In the forward pass (when $L$ is computed), we will use the Heaviside function to obtain $x_k$.
In the backward pass (when computing the gradient $\nabla L$), we will assume that $x_k$ was obtained using a
hot sigmoid with temperature $\tau$ (the 'smooth Heaviside function').
Figure~\ref{sigmoid} shows the shape of this function for different values of temperature $\tau$,
and Figure~\ref{sigmoid-df} shows its derivative.
The closer $\tau$ is to zero, the more closely the sigmoid fits the stepped Heaviside function. 
In fact, the following rule is used:
\begin{equation}
x = H(t),~~~~~~~~~~~~~~\frac{\partial x}{\partial t} = \frac{\partial \sigma(t/\tau)}{\partial t}.
\end{equation}
During the training the temperature is usually assumed to be high $\tau \ge 1$ at first
and the derivatives of discrete $x_k$ by the logits $t_k$
will be different from zero (while $x_k$ remains discrete). 
As the temperature gradually decreases, the derivatives outside the region $t \sim 0$ decrease, but they nevertheless 'steer' 
the values of the discrete parameters $x_k$ to their optimal values.
The software implementation of this computation is discussed in Appendix~\ref{AppC}.

\subsection{The Stochastic Sigmoid}

It is often convenient to consider the binary parameter $x\in\{0,1\}$ as a random variable that takes the value 1 
with probability $p$ and 0 with probability $1-p$.
This nondeterminism helps to get out of local minima, and
in supervised learning tasks leads to data augmentation, thus reducing overfitting.

To construct a stochastic estimator of the discrete value $x$ from the real logit $t$,
let us use the logistic probability distribution $\mathcal{L}(0,s)$:
\begin{equation}
f(t) = \frac{e^{-t/s}}{(1+e^{-t/s})^2},~~~~~~~~~~~~~F(t) = \int\limits^t_{-\infty}f(z)\, dz = \sigma(t/s),
\end{equation}
where $f(t)$ is the probability density function, $F(t)$ is the probability distribution, and $s$ 
is a parameter characterizing the dispersion of the random variable.

The probability that the random variable $\varepsilon \sim \mathcal{L}(0,s)$ is less than given $t$,
is by definition equal to $P(\varepsilon < t) = F(t) = \sigma(t/s)$. 
This means that with probability $\sigma(t/s)$ the value of $t-\varepsilon$ is positive 
(the probability density $f(t)$ is symmetric and $\varepsilon$ can be both subtracted and added to $t$).
Therefore, if one subtracts $\varepsilon\sim \mathcal{L}(0,s)$ from the logit in a deterministic sigmoid with temperature 
and computes the Heaviside function from the result, then with probability $\sigma(t/s)$ one will obtain 1, and 0 otherwise:
\begin{equation}
\varepsilon\sim \mathcal{L}(0,s)~~~\Rightarrow~~~H(t-\varepsilon) = 
\left\{ 
\begin{array}{lll}
1 & \text{with probability} & \sigma(t/s)\\
0 & \text{otherwise}
\end{array}
\right.
\end{equation}
Note also that the logistic distribution 
is a limiting case of the Gumbel distribution~(\cite{maddison2017concrete}, \cite{jang2017categorical}) when 
dealing with a categorical variable that has only two values.

The value of the parameter $s$ depends on the problem we are trying to solve. We can set $s=1$;
then for the temperature annealing (decreasing $\tau$) the stochasticity will decrease only 
if the absolute values of logits grow large during the training.
For some problems, the 'stochasticity with heating' is more appropriate,
where $s$ is minimal at the beginning of training, then increases to a maximum value and decreases towards the end of training.

As a result, we have two 'temperatures' with different annealing schedules. 
One temperature ($\tau$) regulates the transition from continuous variable space to discrete, 
and the second one ($s$) regulates the degree of random variability of the parameters.

\subsection{CONGA - CONstrained Gradient descent with Adaptation}

The optimization problem in the presence of constraints leads to the loss function (\ref{loss}) with a cusp 
(the penalty is zero when $w(\mathbf{x})\ge 0$, and different from zero otherwise).  
One problem in applying the gradient method to such loss function is related to the choice of the parameter $\gamma$. 
If it is small, the penalty may be insufficient 
and the minimum point will be located in the forbidden region $w(\mathbf{x}) > 0$. 
When $\gamma$ is large, the solution on the boundary of the region starts to experience oscillations, which slows down convergence. 
The problem is partially solved if we apply the exponential moving average (EMA) smoothing to the gradients.
Nevertheless, it would be advantageous to have a method that selects the optimal value of the hyperparameter $\gamma$.
In what follows, we will put forward such a method. 

Following~\cite{li2021novel} we require that in the next step of the gradient method
the positive penalty for violating the constraint decreases: $w(\mathbf{x}-\lambda\,\nabla L) = (1-\mu) \,w(\mathbf{x})$.
This leads to a new adaptive method to select the parameter $\gamma$, which we have named
CONGA (CONstrained Gradient descent with Adaptation):
\begin{equation}\label{gamma-CONGA}
\gamma  = \frac{\nabla v \nabla w+\mu \,w}{w^{\nu-1}\, (\nabla w)^2},~~~~~~\gamma=\max(0,\,\gamma),
\end{equation}

In fact, instead of the hyperparameter $\gamma$, another one $\mu\in [0...1]$ is introduced.
It is worth pointing out that unlike $\gamma$, this new parameter has the simple meaning and the corresponding rule to select 
its value: the target relative 
reduction of the violation of the constraint condition $w(\mathbf{x})$
at the next step of the gradient method.
Additionally, we introduce the parameters $\beta_v,~\beta_w\in[0...1]$ of EMA smoothing of 
gradients $\nabla v$ and $\nabla w$ using them in (\ref{gamma-CONGA}) and in the loss-function gradient $\nabla L$.
Details of the computations and comparisons with the adaptive  method of~\cite{li2021novel} are given
in the Appendix~\ref{AppCONGA}.

\subsection{The Population of Solutions}

Using the GPU one can efficiently calculate the functions $v(\mathbf{x})$, $w(\mathbf{x})$ and their gradients
simultaneously for several trajectories in the parameter space.
In our experiments, we create $m$ 'individuals' with logits $\mathbf{t}^{(k)}$, $k=1,...,m$,
which have different initial positions.
These vectors are independently and simultaneously varied according to the CONGA method.

In addition to the 'innate' (initial) diversity, individuals have different values of the hyperparameter $\mu$. As experiments have shown, the hyperparameter $\mu$ for each example has its own value at which the given problem better converges to the optimal solution. Therefore, one of the tasks of a population of 'individuals' is to choose the optimal distribution of a given hyperparameter.
The lifetime of a population of individuals is limited by another hyperparameter (\verb|epochs|), after this interval the selection of 'individuals' is performed in the following way.
The best 'individuals' (top 20\%) are selected that have the largest maximum backpack value. The range of parameter $\mu_{best}\in[\mu_{min},\mu_{max}]$ is determined for these solutions, and then the range boundaries are proportionally expanded: 
$\mu_{next}\in[\mu_{min}/frac,\mu_{max}\cdot frac]$. The resulting range is used to create new 'individuals' of the population.
As the result, the directed evolutionary dynamics is obtained.

\begin{algorithm}
\caption{CONGA}\label{alg:CONGA}
\begin{algorithmic}
\State $gen\_id = 0$
\ForAll{n\_generations}
\If{$gen\_id > 0$} 
\State $p_k = selection(p_{k-1})$
\Comment{selecting top individuals from generation}
\Else \State $p_k \gets$ individual
\EndIf
\ForAll{epochs}
\ForAll{individual}
\State $x_k = hot\_sigmoid(p_k)$
\Comment{soft/hard sigmoid}
\State $v_k = v(\mathbf{x}_k)$     
\State $w_k = w(\mathbf{x}_k)$     
\State $\mathbf{V}_k ~\gets \beta_v\,\mathbf{V}_k ~\,+ (1-\beta_v) \nabla v(\mathbf{x}_k)$
\Comment{gradient smoothing}
\State $\mathbf{W}_k \gets \beta_w\,\mathbf{W}_k + (1-\beta_w) \nabla w(\mathbf{x}_k)$   
\State $\gamma_k  = \max(0,~(\mathbf{V}_k\mathbf{W}_k+\mu \,w_k/\lambda ) / (w^{\nu-1}_k\, \mathbf{W}^2_k)\,)$
\State $\mathbf{p}_{k+1} \gets \mathbf{p}_k - \lambda \,(-\mathbf{V}_k + \gamma_k\,\max(0, w_k)^{\nu-1}\,W_k)$         
\EndFor
\EndFor
\EndFor
\end{algorithmic}
\end{algorithm}

\section{Experimental Results}\label{experiments}

We will demonstrate the performance of the CONGA method using the dataset for 0--1 Knapsack Problem described in~\cite{instances01kp}. The dataset is split into groups depending on the number of items and the correlation of values and their weights as shown in the table~\ref{table:datasets}. The following hardware platform was used in the experiments: Intel(R) Xeon(R) 2.00GHz CPU, Tesla T4 GPU.
The source code is implemented in Python and is available online in the repository 
at \href{https://github.com/QuDataAI/CONGA}{https://github.com/QuDataAI/CONGA}. 

\begin{table}
\centering
\caption{Datasets.}
\label{table:datasets}
\begin{tblr}{
  colspec = {l|l|l|l}  
}
\hline[1.2pt]
name & items & correlation & tasks\\
\hline[1.2pt]
LD-UC & 4-20 & uncorrelated & 10\\ 
\hline[1pt]
HD-UC & 100-10000 & uncorrelated & 7\\ 
\hline[1pt]
HD-WC & 100-10000 & weakly correlated & 7\\ 
\hline[1pt]
HD-SC & 100-10000 & strongly correlated & 7\\ 
\hline[1.2pt]
\end{tblr}
\end{table}

\subsection{Single-Agent Case}
Let us now consider the application of this method without directed evolution in the mode with only one agent.
For this case, the number of populations and the number of agents are $n\_generations=1, n\_agents=1$.

The agent, having received the task, initializes the initial state of the backpack packing vector randomly with a normal distribution. Then, using the CONGA method, the value of the backpack is iteratively optimized while changing the packing vector and taking into account the weight constraint. The maximum number of such iterations is limited by the hyperparameter $epochs=2000$ in order to prevent an infinite search for the optimal solution.

As a result of optimization, the agent can come to either an exact or an approximate solution, however the iteration always 
converges to a valid solution.

To evaluate the results we will use the following metrics:
\begin{equation}
acc = \frac{1}{n}\sum_{i=1}^{n} (v_{opt} == v),~~~~~time = \frac{1}{n}\sum_{i=1}^{n} t
\end{equation}

where $v_{opt}$ is the optimal value of the backpack, $v$ is the value of the backpack after optimization, $n$ is the number of tasks in the dataset, $t$ is the time to solve one problem in seconds. Thus, $acc$ shows how many problems reached the optimal solution, and $time$ is the average time to solve one problem in seconds.

By running this experiment, we obtained the following metrics: $\{acc:0.48,~~time:2.655\}$, therefore, we solve optimally  
less than half of the problems, while spending approximately 3 seconds on each problem.

\begin{figure}
\centering  
\begin{tikzpicture}
\begin{axis}  
[  
    ybar,   
    enlargelimits=0.12,  
    legend style={at={(0.3,-0.25)},  
      anchor=north,legend columns=-1}, 
    ylabel={acc}, 
    symbolic x coords={LD-UC, HD-UC, HD-WC, HD-SC},  
    xtick=data,  
    nodes near coords,  
    nodes near coords align={vertical},  
]  
\addplot coordinates {(LD-UC, 0.80) (HD-UC, 0.14) (HD-WC, 0.29) (HD-SC, 0.29)}; 
\end{axis}  
\end{tikzpicture}  
\caption{Single-agent accuracy for different task types}\label{fig:single}
\end{figure}
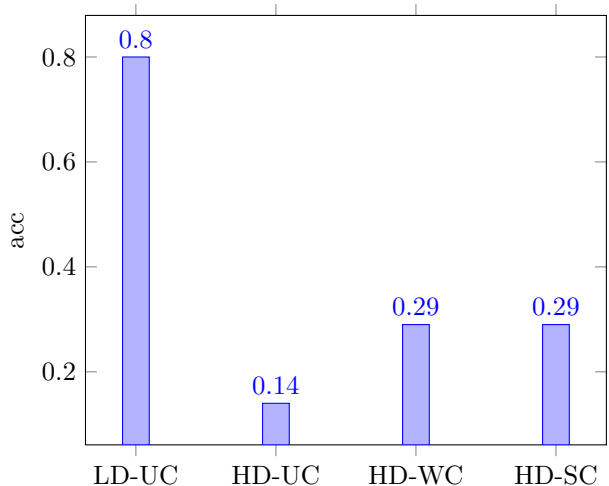

The outcome of the optimizer application to different types of problems is presented in Fig.~\ref{fig:single}.
As one can see, the method works better with a dataset in which there are fewer objects, but one agent is clearly not enough to solve more complex problems.

\subsection{Multi-Agent Case}

The multi-agent approach allows one to solve the problem in parallel by several agents, as representatives of the population that has a common goal, but each agent has a unique initial state of the backpack packing vector. Thanks to GPU parallelization, we can significantly increase 
the number of agents without losing processing speed.

\begin{figure}   
\centering
\begin{tikzpicture}
  \begin{axis}[ 
    xlabel=$agents$,
    ylabel={$acc$}
  ] 
    \addplot [
    color=blue,
    mark=square,
    ]
    coordinates {(1,0.4838)(10,0.6774)(25,0.7419)(50,0.7742)(100,0.7419)};
  \end{axis}
\end{tikzpicture}
\caption{Metric $acc$ dependence on the number of agents}\label{fig:Multi-acc}
\end{figure}
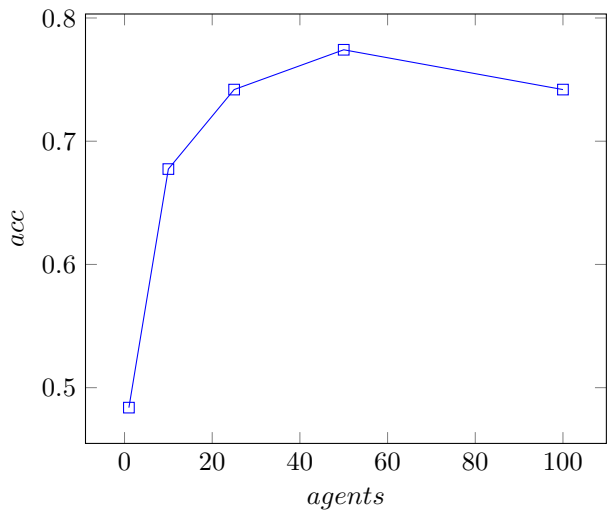

We can see how the optimizer's result changes depending on the number of agents in Fig.~\ref{fig:Multi-acc}.

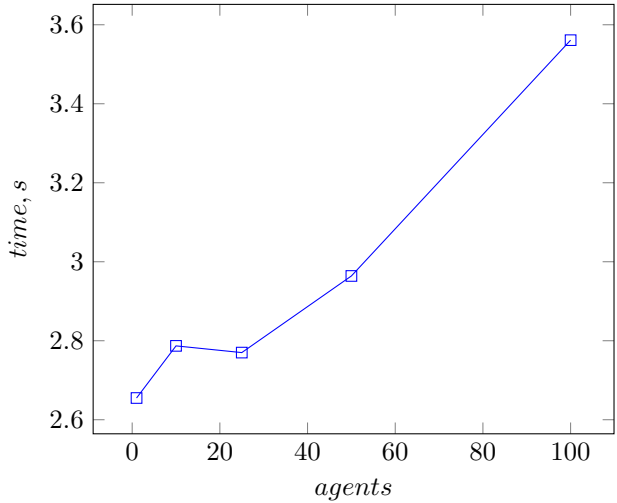
\begin{figure} 
\centering
\begin{tikzpicture}
  \begin{axis}[ 
    xlabel=$agents$,
    ylabel={$time,s$}
  ] 
    \addplot [
    color=blue,
    mark=square,
    ]
    coordinates {(1,2.655)(10,2.787)(25,2.770)(50,2.964)(100,3.561)};
  \end{axis}
\end{tikzpicture}
\caption{Dependence of the mean solution time on the number of agents} \label{fig:Multi-time}  
\end{figure} 

As one can see, 50 agents is optimal for solving problems, and yields the metric values equal to $\{acc:0.74(+0.26),~~time:2.964(+0.309)\}$, while the average time spent on solving one problem increased by only $309 ms$ thanks to parallel computing on the GPU.

A further increase in the number of agents linearly increases the time to process one task as shown in Fig.~\ref{fig:Multi-time}, 
but the result does not improve. Thus, we fix the hyperparameter $n\_agents=50$

\subsection{Directed Evolutionary Dynamics}

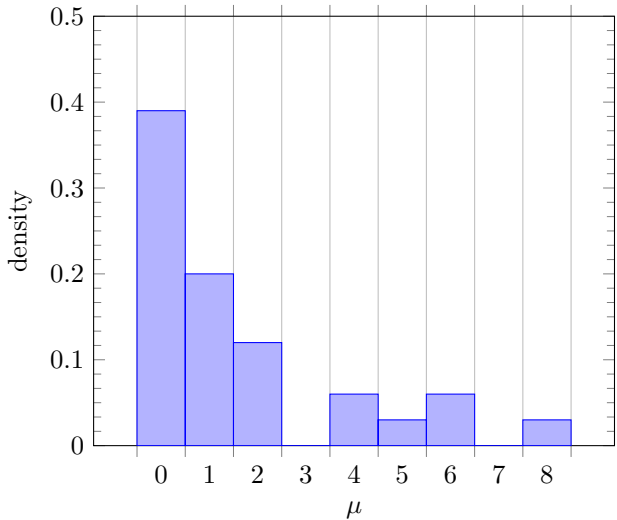
\begin{figure}
\centering
\begin{tikzpicture}
\begin{axis}[xlabel={$\mu$},ylabel={density},
ybar interval, ymax=0.5,ymin=0, minor y tick num = 5]
\addplot coordinates { (0, 0.39) (1, 0.2) (2, 0.12) (3, 0) (4, 0.06) (5,0.03) (6,0.06) (7,0) (8,0.03) ( 9,0.0)};
\end{axis}
\end{tikzpicture}
\caption{Histogram of parameter density $\mu$ among optimal solutions} \label{fig:Multi-hist}
\end{figure}

Analyzing the process of solving problems in the dataset, we noticed that some problems are successfully solved with small values of the parameter $\mu$ from Eq.~\ref{gamma-CONGA}, while others require higher values as shown in the histogram~\ref{fig:Multi-hist}.

Therefore, we proceeded to endow our method with the properties of a genetic algorithm. 
Agents from the population randomly select a hyperparameter value $\mu$ from the range $[0.2,8.0]$. 
After reaching the maximum number of epochs, agents with the best backpack values are selected. 
As a result of selection, we determine a new hyperparameter range $\mu$ that better optimizes this problem. 
Then we create a new population of agents that inherit the updated range $\mu$ and therefore optimize the solution based on the experience of their ancestors. The number of such iterations is determined by the hyperparameter $n\_generations$.

Executing such a process of evolutionary directed dynamics, the CONGA method finds the optimal solution for all problems in the dataset;
correspondingly, we obtain the following metrics: $\{acc:1.0(+0.26),~~time:6.132(+3.168)\}$.
 
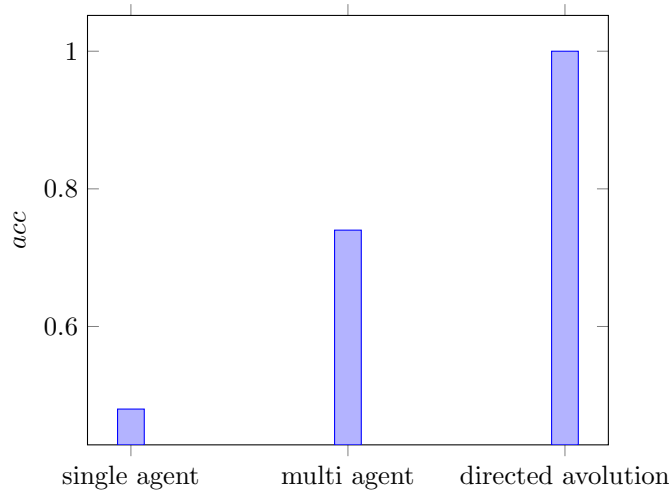
\begin{figure}
\centering
\begin{tikzpicture}
    \begin{axis}[
    		ylabel={$acc$},
            ybar,
            symbolic x coords={single agent,multi agent,directed avolution},
            xtick=data,
        ]
        \addplot 
coordinates { (single agent, 0.48) (multi agent, 0.74) (directed avolution, 1.0)};
    \end{axis}
\end{tikzpicture}
\caption{Dependence of the result on the given experiments}\label{fig:total-res}
\end{figure}

Thus, in the Fig.~\ref{fig:total-res} we see the contribution of each experiment to the final result.

The detailed results of this experiment in comparison with other approximate methods are presented in the Appendix~\ref{AppResults}.

\section{Conclusions}\label{conclusion}

In conclusion, this paper addressed the challenge of optimizing discrete parameters within defined constraints. 
Introducing the stochastic sigmoid with temperature, we proposed the novel adaptive gradient method, CONGA, employing 
a population-based dynamical evolution approach. 
Each individual within the population undergoes variation based on environmental gradients, characterized by two 
temperature parameters with distinct annealing schedules. 
The evolutionary dynamics involve the elimination of unadapted individuals, while optimal ones interbreed, fostering a 
directed evolution. 
While in general case the method finds approximate solutions, our experiments have shown that in many cases 
the method yields optimal solutions, which makes it one step closer to exact methods. 
Evolutionary dynamics adapts its parameters to a specific problem, allowing one to find rather non-trivial solutions.
The efficacy of the proposed method was demonstrated through its application to the classic combinatorial problem of 
optimal packing of a backpack (0–1 Knapsack Problem).


\section{Acknowledgments}

We wish to thank the authors of papers~\cite{li2021novel,mena2018learning} for helpful communication on their methods and datasets.
A special thanks to Gloria Estefan for the hit song {\em Conga}, the unbounded source of energy that kept our algorithms running.

\printbibliography

\newpage
\appendix

\section{Code Implementation for Hot Sigmoid}\label{AppC}

\begin{wrapfigure}{r}{0.5\textwidth}\centering
    \includegraphics[width=0.45\textwidth]{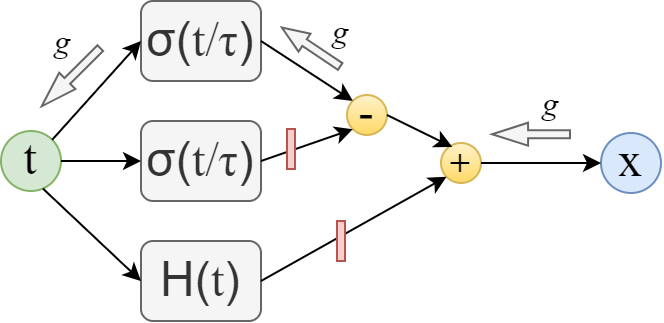}     
    \caption{\small Computation graph for the deterministic hot sigmoid}\label{graph}
\end{wrapfigure}

We present here the program code for computing the stochastic sigmoid with temperature.
Instead of the Heaviside function we will use rounding of the sigmoid output, which leads to the same result 
(and it is faster in the pytorch framework).
The values of \verb|x_hard| are, while real, but nevertheless discrete values of $0$ or $1$.
When calculating this quantity, there is a detachment from the graph (\verb|detach|), so the gradient will not propagate along this branch.
Then we add the difference between the numerically equal values of \verb|x_soft| and \verb|x_soft.detach()|, which does not change the result
(the grouping by parentheses around the difference is important to avoid the loss of precision).
Since one of the terms (\verb|x_soft|) is not detached from the graph, the gradient computed from the sigmoid with temperature will follow it.
In  Figure~\ref{graph}, detachment from the graph is denoted by the red vertical bars. 
The forward computation runs left-to-right along the graph, and the gradient $g$ propagates right-to-left.

{\large
\begin{lstlisting}[language=Python]
def hot_sigmoid(t, tau=1., s=1., hard=True, rand=True, eps=1e-8):
    if rand:        
        r = torch.rand_like(t)                   # ~ U(0,1)
        r = torch.log(eps + r / (1-r+eps))       # ~ L(0,1)
        t = t + s * r                            # + L(0,s)
 
    x_soft = torch.sigmoid(t/tau)
    if not hard:
        return x_soft
     
    x_hard = x_soft.detach().round()
    return x_hard + (x_soft - x_soft.detach())
\end{lstlisting}
}

If the \verb|hard| parameter is \verb|True|, the function will return discrete $x$, and calculate the gradient using a smoothed 
version of the step-function. In the \verb|hard=False| mode, the sigmoid will also be used to compute $x$ directly.
For generality, we will also combine the deterministic and stochastic functions into one (parameter \verb|rand|).

\newpage
\section{CONstrained Gradient descent with Adaptation}\label{AppCONGA}

Suppose that the change of parameters leads to escape into the forbidden region of the search space (a penalty is triggered): 
$w(\mathbf{x})> 0$.
Then the loss function (\ref{loss}) and its gradient would have the form:
\begin{equation}
L = -v + \frac{\gamma}{\nu}\,w^\nu,~~~~~~~\nabla L = -\nabla v + \gamma\,w^{\nu-1}\, \nabla w.
\end{equation}
We choose $\gamma$ so that $w(\mathbf{x})$ decreases in the next step of the gradient method:
\begin{equation}
w(\mathbf{x}-\lambda\,\nabla L) = (1-\mu) \,w(\mathbf{x}),
\end{equation}
where $\mu \in [0...1]$ is a constant (another hyperparameter). The closer $\mu$ is to unity, the more strongly $w$ decreases.
Assuming that the change in the parameters $\lambda\,\nabla L$ is small, let us expand the left part of the equality into a Taylor series: 
$w(\mathbf{x}-\lambda\,\nabla L) \approx w(\mathbf{x}) -\lambda\,\nabla L\,\nabla w(\mathbf{x})$
and substitute $\nabla L$. This yields the following value for the factor at penalty:
\begin{equation}
\gamma  = \frac{\nabla v \nabla w+\mu \,w}{w^{\nu-1}\, (\nabla w)^2},~~~~~~\gamma=\max(0,\,\gamma).
\end{equation}
Trimming the value of $\gamma$ from below is necessary to ensure that the weight of the penalty is not negative when $\nabla v \nabla w < 0$.
The possible singularity in the denominator is eliminated, as usual, by adding a small $\text{eps}\sim 10^{-6}$.

In the method described in~\cite{li2021novel}, it is required that along with the decrease of $w$, the target function $v$ increases.
The corresponding AGA (Adaptive Gradient Ascent) algorithm is as follows:
\begin{equation}\label{AGA}
\gamma = \left\{
\begin{array}{llll}
 0             & \nabla v \nabla w < 0                  & (1)  \\
 (r_1+r_2)/2   & \nabla v \nabla w > 0 ~~\&~~ r_1 < r_2 & (2)  \\
 r_3           & \text{otherwise}                       & (3)  
\end{array}
\right.
\end{equation}
where 
$$
r_1 = \frac{ \nabla v \nabla w}{w^{\nu-1}\,(\nabla w)^2 },~~~~~
r_2 = \frac{(\nabla v)^2}{ w^{\nu-1}\, \nabla v \nabla w},~~~~~
r_3 = \frac{\nabla v \nabla w + w/\lambda }{w^{\nu-1}\, (\nabla w )^2}.
$$
In fact, when $\nabla v \nabla w > 0$, it is almost always $r_1 < r_2$, since this inequality is not satisfied 
only when the gradients of $\nabla v$ and $\nabla w$ are exactly parallel (or antiparallel).
Therefore, the third condition is rarely satisfied.

\begin{figure}[!ht]
    \includegraphics[width=15cm]{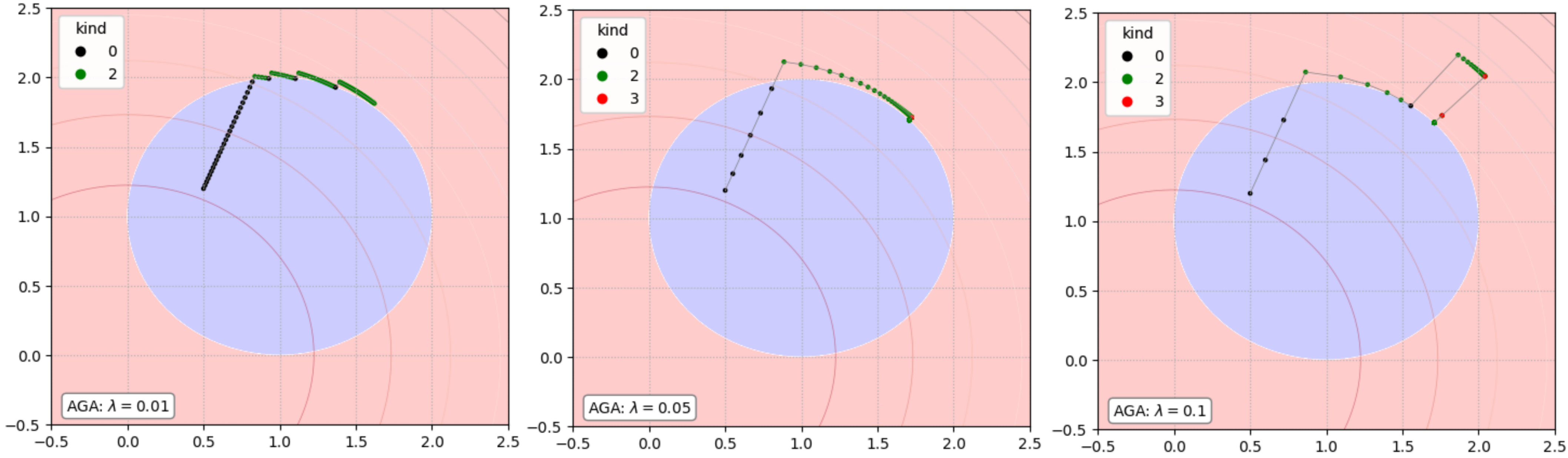}     
    \caption{\small 
    Search for the optimal solution by the AGA method for different $\lambda$. 
    Throughout the following, the target function $v = x^2+y^2= \max$ and $\mathbf{x}=\{x,y\}$ are assumed to be continuous. 
    The constraint function:  $w = (x-1)^2 + (y-1)^2 - 1 \le 0$ }\label{aga_01}
\end{figure}

Note that the original AGA method is quite sensitive to the learning rate $\lambda$.
Figure~\ref{aga_01} shows the trajectories of the gradient method from the same point at different $\lambda$.
The blue circle outlines the region covering the allowed solutions. The color of the points corresponds to the three conditions 
of the AGA method in (\ref{AGA}) and the zero index is given to the points inside the allowed region.

\begin{figure}[!ht]    
    \includegraphics[width=15cm]{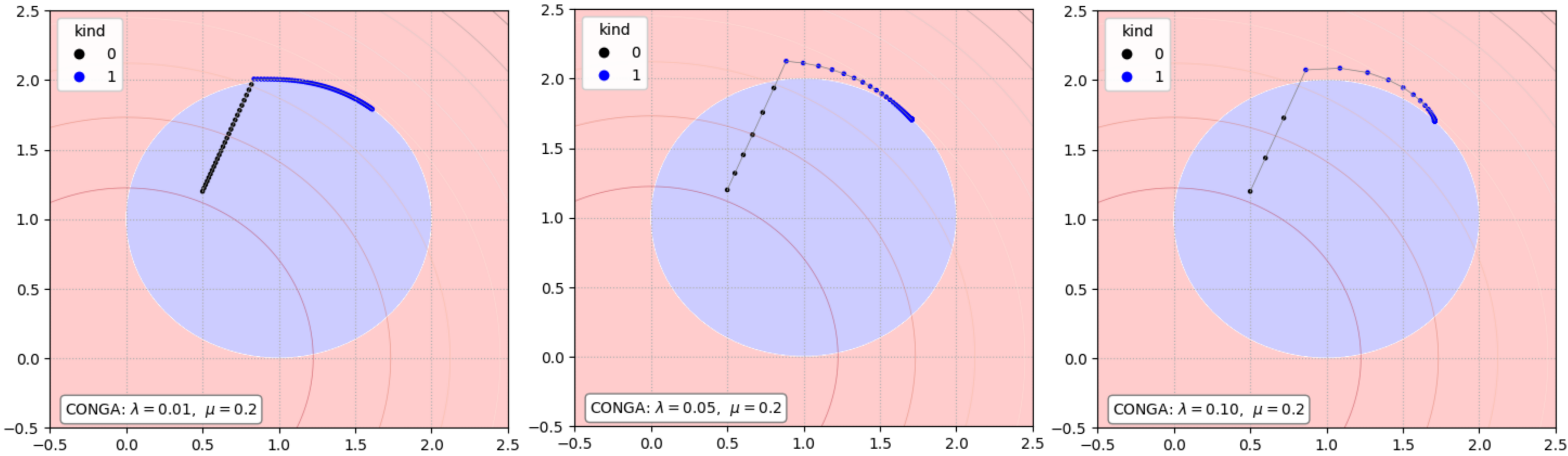}     
    \caption{\small Search for the optimal solution by the CONGA method for different $\lambda$.}\label{conga_01}
\end{figure}

The CONGA method proposed in this paper is significantly less sensitive to the choice of learning rate (see Fig.~\ref{conga_01}).

\begin{figure}[!ht]    
    \includegraphics[width=15cm]{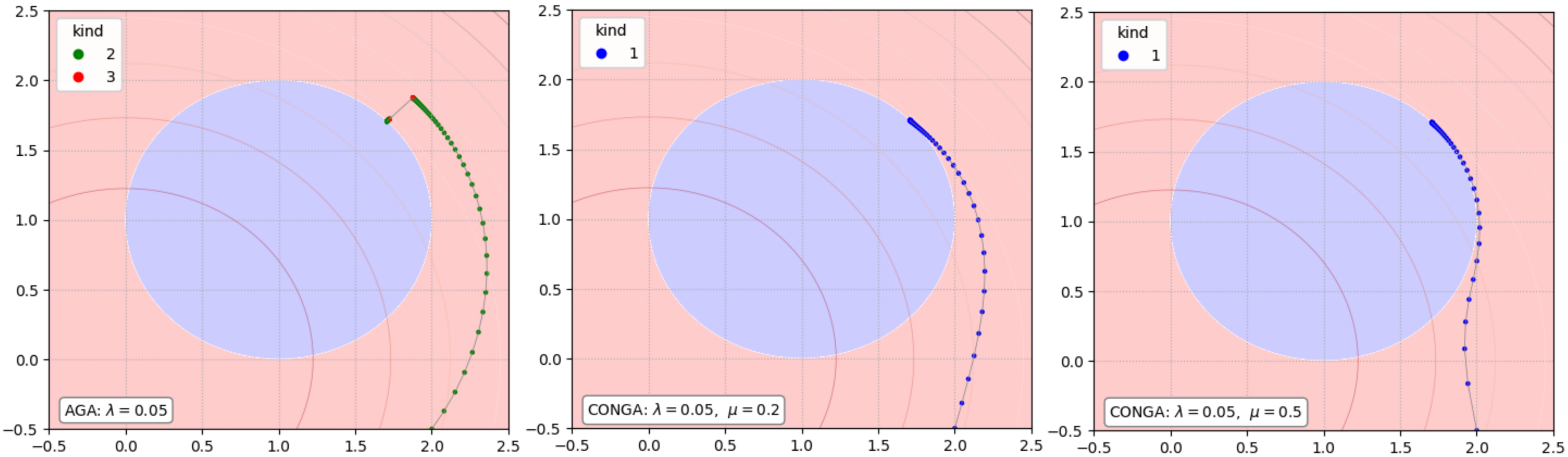}     
    \caption{\small Search for the optimal solution starting from the forbidden region: AGA method compared to CONGA (for 2 
    different values of $\mu$). The closer $\mu$ is to $1$, the more strongly the solution is pushed to the boundary.}\label{conga_02}
\end{figure}

The second advantage of CONGA comes from the presence of a hyperparameter $\mu$, which regulates the desired 'attraction'
of the iterative procedure to the interior of the allowed region.
Figure~\ref{conga_02} shows a solution that starts in the forbidden region. 
The first graph corresponds to the AGA method, and the second and third to different values of the hyperparameter $\mu$ in the CONGA method.

\begin{figure}[!ht]    
    \includegraphics[width=15cm]{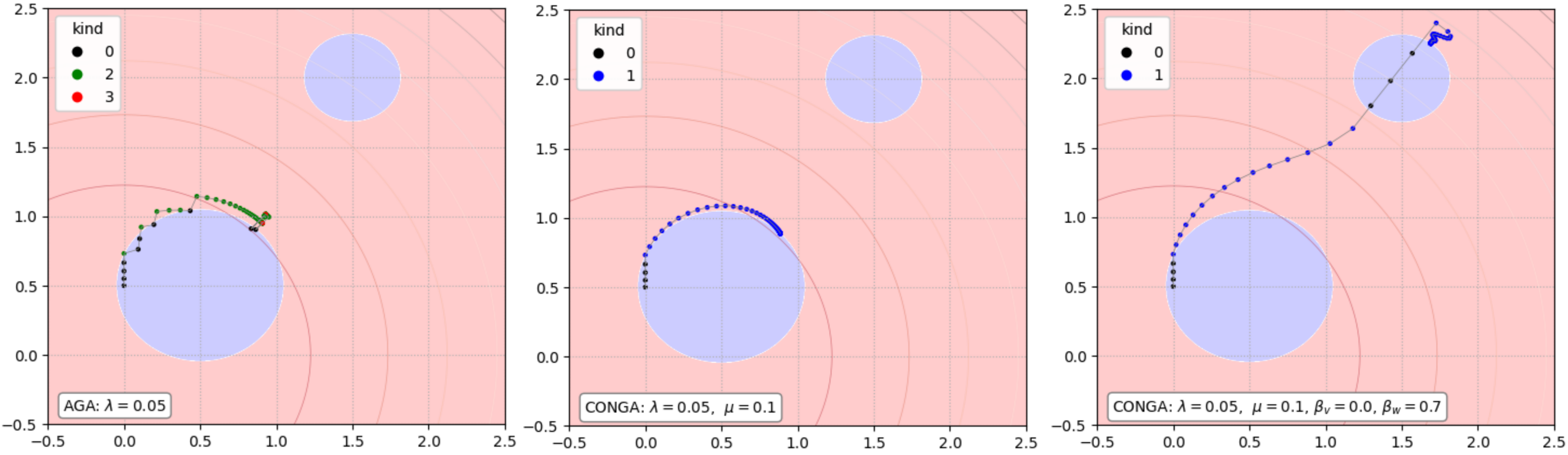}     
    \caption{\small Search for an optimal solution for a non-convex constraint function by AGA, CONGA w/o EMA and CONGA with EMA methods.}
    \label{aga_conga_2constr}
\end{figure}

Another important feature of CONGA is the introduction of exponential smoothing for the gradient $\nabla v$ 
(with the parameter $\beta_v\in[0...1]$)
and $\nabla w$ (with the parameter $\beta_w\in[0...1]$): 
\begin{equation}
V_{k+1} = \beta_v \, V_k + (1-\beta_v)\, \nabla v_k,~~~~~~~~~~W_{k+1} = \beta_w \, W_k + (1-\beta_w)\, \nabla w_k
\end{equation}
The closer $\beta_w$ is to $1$, the deeper the iteration can penetrate into the forbidden region, while still eventually 
converging to the valid solution.
This allows us to find an optimal solution even in situations where the constraint is not a convex function (Figure~\ref{aga_conga_2constr}).

\newpage
\section{Summary of CONGA Method Application}\label{AppResults}

We present here the comparison of our results to several other algorithms. To compare the results of the CONGA method proposed in this paper,
the following algorithms for solving the 0--1 KP problem were used: Branch and Bound (BB), Genetic Algorithm (GA), Greedy Search Algorithm (GSA) and Simulated Annealing (SA).
One recent review of these algorithms application to the dataset we use to test our method is given in~\cite{Ezugwu2019}.
Therefore, in this work, we solved the problems of the above dataset using the CONGA algorithm and compared it with benchmark results collected in~\cite{Ezugwu2019}.

The results of comparing the CONGA method with other algorithms by type of problem are presented in the Tables \ref{table:result-ld-uc} (LD-UC datasets), \ref{table:result-hd-uc} (HD-UC datasets), \ref{table:result-hd-wc} (HD-WC datasets), and \ref{table:result-hd-sc} (HD-SC datasets).
The values of the hyperparameters used in our experiments are given in Table~\ref{table:hyperparameters}.

\begin{table}[hbt]
\centering
\caption{Hyperparameters}\label{table:hyperparameters}
\begin{tblr}{
  colspec = {ll}  
}
\hline[1.2pt]
Hyper-parameter & Values & Description\\
\hline[1pt]
n\_generations & 2 & number of generations\\ 
epochs & 2000 & maximum epochs per generation\\ 
n\_agents & 50 & number of agents in generation \\ 
lr & 1e-1 & learning rate \\
nu & 1.0 & fraction of boundary part \\
mu1 & 0.2 & minimum value of mu for initial distribution \\
mu2 & 0.8 & maximum value of mu for initial distribution \\
beta\_v & 0.5 & EMA beta for values gradient \\
beta\_w & 0.5 & EMA beta for weights gradient \\
eps & 1e-6 & avoid zero div in gamma calc \\
tau1 & 30 & initial temperature for hot sigmoid \\
tau\_hot & 30 & hot temperature value at tau\_warmap\_epochs  \\
tau2 & 0.01 & final temperature value at tau\_max\_epochs \\
tau\_warmap\_epochs  & 1  & epochs where temperature reaches tau\_hot \\
tau\_max\_epochs & 2000 & epochs where temperature reaches tau2 \\
\hline[1.2pt]
\end{tblr}
\end{table}

\begin{table}[!htbp]
\centering
\caption{Results obtained on LD-UC datasets.}\label{table:result-ld-uc}
\begin{tblr}{
  colspec = {llcccc}  
}
\hline[2pt]
  Algorithm &   Dataset &  Items & Best Value & Optimal Value & Shortfall \\
\hline[2pt]
BB & f1\_1-d\_kp\_10\_269 & 10 & 280 & 295 & 15 \\
\textbf{GA} & f1\_1-d\_kp\_10\_269 & 10 & 295 & 295 & \textbf{0} \\
GSA & f1\_1-d\_kp\_10\_269 & 10 & 288 & 295 & 7 \\
SA & f1\_1-d\_kp\_10\_269 & 10 & 294 & 295 & 1 \\
\textbf{CONGA} & f1\_1-d\_kp\_10\_269 & 10 & 295 & 295 & \textbf{0} \\
\hline[1pt]
BB & f2\_1-d\_kp\_20\_878 & 20 & 972 & 1024 & 52 \\
\textbf{GA} & f2\_1-d\_kp\_20\_878 & 20 & 1024 & 1024 & \textbf{0} \\
\textbf{GSA} & f2\_1-d\_kp\_20\_878 & 20 & 1024 & 1024 & \textbf{0} \\
SA & f2\_1-d\_kp\_20\_878 & 20 & 972 & 1024 & 52 \\
\textbf{CONGA} & f2\_1-d\_kp\_20\_878 & 20 & 1024 & 1024 & \textbf{0} \\
\hline[1pt]
BB & f3\_1-d\_kp\_4\_20 & 4 & 24 & 35 & 11 \\
\textbf{GA} & f3\_1-d\_kp\_4\_20 & 4 & 35 & 35 & \textbf{0} \\
GSA & f3\_1-d\_kp\_4\_20 & 4 & 28 & 35 & 7 \\
\textbf{SA} & f3\_1-d\_kp\_4\_20 & 4 & 35 & 35 & \textbf{0} \\
\textbf{CONGA} & f3\_1-d\_kp\_4\_20 & 4 & 35 & 35 & \textbf{0} \\
\hline[1pt]
BB & f4\_1-d\_kp\_4\_11 & 4 & 22 & 23 & 1 \\
\textbf{GA} & f4\_1-d\_kp\_4\_11 & 4 & 23 & 23 & \textbf{0} \\
GSA & f4\_1-d\_kp\_4\_11 & 4 & 19 & 23 & 4 \\
SA & f4\_1-d\_kp\_4\_11 & 4 & 22 & 23 & 1 \\
\textbf{CONGA} & f4\_1-d\_kp\_4\_11 & 4 & 23 & 23 & \textbf{0} \\
\hline[1pt]
BB & f5\_1-d\_kp\_15\_375 & 15 & 469,16 & 481,07 & 11,91 \\
GA & f5\_1-d\_kp\_15\_375 & 15 & 477 & 481,07 & 4,07 \\
\textbf{GSA} & f5\_1-d\_kp\_15\_375 & 15 & 481,07 & 481,07 & \textbf{0} \\
\textbf{SA} & f5\_1-d\_kp\_15\_375 & 15 & 481,07 & 481,07 & \textbf{0} \\
\textbf{CONGA} & f5\_1-d\_kp\_15\_375 & 15 & 481,07 & 481,07 & \textbf{0} \\
\hline[1pt]
BB & f6\_1-d\_kp\_10\_60 & 10 & 49 & 52 & 3 \\
\textbf{GA} & f6\_1-d\_kp\_10\_60 & 10 & 52 & 52 & \textbf{0} \\
GSA & f6\_l-d\_kp\_10\_60 & 10 & 43 & 52 & 9 \\
\textbf{SA} & f6\_l-d\_kp\_10\_60 & 10 & 52 & 52 & \textbf{0} \\
\textbf{CONGA} & f6\_l-d\_kp\_10\_60 & 10 & 52 & 52 & \textbf{0} \\
\hline[1pt]
BB & f7\_1-d\_kp\_7\_50 & 7 & 96 & 107 & 11 \\
\textbf{GA} & f7\_1-d\_kp\_7\_50 & 7 & 107 & 107 & \textbf{0} \\
GSA & f7\_l-d\_kp\_7\_50 & 7 & 100 & 107 & 7 \\
SA & f7\_1-d\_kp\_7\_50 & 7 & 102 & 107 & 5 \\
\textbf{CONGA} & f7\_1-d\_kp\_7\_50 & 7 & 107 & 107 & \textbf{0} \\
\hline[2pt]
\end{tblr}
\end{table}

\begin{table}[!htbp]
\centering
\caption{Results obtained on HD-UC datasets.}\label{table:result-hd-uc}
\begin{tblr}{
  colspec = {llcccc}  
}
\hline[2pt]
  Algorithm &   Dataset &  Items & Best Value & Optimal Value & Shortfall \\
\hline[2pt]
\textbf{SA} & knapPI\_1\_100\_1000\_1 & 100 & 9147 & 9147 & \textbf{0} \\
\textbf{GA} & knapPI\_1\_100\_1000\_1 & 100 & 9147 & 9147 & \textbf{0} \\
GSA & knapPI\_1\_100\_1000\_1 & 100 & 2983 & 9147 & 6164 \\
BB & knapPI\_1\_100\_1000\_1 & 100 & 8026 & 9147 & 1121 \\
\textbf{CONGA} & knapPI\_1\_100\_1000\_1 & 100 & 9147 & 9147 & \textbf{0} \\
\hline[1pt]
SA & knapPI 1 200\_1000\_1 & 200 & 10163 & 11238 & 1075 \\
GA & knapPI 1 200\_1000\_1 & 200 & 102340 & 11238 & -91102 \\
GSA & knapPI 1 200\_1000\_1 & 200 & 4544 & 11238 & 6694 \\
BB & knapPI\_1\_200\_1000\_1 & 200 & 10436 & 11238 & 802 \\
\textbf{CONGA} & knapPI\_1\_200\_1000\_1 & 200 & 11238 & 11238 & \textbf{0} \\
\hline[1pt]
SA & knapPI\_1\_500\_1000\_1 & 500 & 21390 & 28857 & 7467 \\
GA & knapPI\_1\_500\_1000\_1 & 500 & 102340 & 28857 & -73483 \\
GSA & knapPI\_1\_500\_1000\_1 & 500 & 9865 & 28857 & 18992 \\
BB & knapPI\_1\_500\_1000\_1 & 500 & 28043 & 28857 & 814 \\
\textbf{CONGA} & knapPI\_1\_500\_1000\_1 & 500 & 28857 & 28857 & \textbf{0} \\
\hline[1pt]
SA & knapPI\_1\_1000\_1000\_1 & 1000 & 36719 & 54503 & 17784 \\
GA & knapPI\_1\_1000\_1000\_1 & 1000 & 130 & 54503 & 54373 \\
GSA & knapPI\_1\_1000\_1000\_1 & 1000 & 14927 & 54503 & 39576 \\
BB & knapPI\_1\_1000\_1000\_1 & 1000 & 53397 & 54503 & 1106 \\
\textbf{CONGA} & knapPI\_1\_1000\_1000\_1 & 1000 & 54503 & 54503 & \textbf{0} \\
\hline[1pt]
SA & knapPI\_1\_2000\_1000\_1 & 2000 & 65793 & 110625 & 44832 \\
GA & knapPI\_1\_2000\_1000\_1 & 2000 & 102340 & 110625 & 8285 \\
GSA & knapPI\_1\_2000\_1000\_1 & 2000 & 25579 & 110625 & 85046 \\
BB & knapPI\_1\_2000\_1000\_1 & 2000 & 109679 & 110625 & 946 \\
\textbf{CONGA} & knapPI\_1\_2000\_1000\_1 & 2000 & 110625 & 110625 & \textbf{0} \\
\hline[1pt]
SA & knapPI\_1\_5000\_1000\_1 & 5000 & 150731 & 276457 & 125726 \\
GA & knapPI\_1\_5000\_1000\_1 & 5000 & 102340 & 276457 & 174117 \\
GSA & knapPI\_1\_5000\_1000\_1 & 5000 & 49306 & 276457 & 227151 \\
BB & knapPI\_1\_5000\_1000\_1 & 5000 & 275720 & 276457 & 737 \\
\textbf{CONGA} & knapPI\_1\_5000\_1000\_1 & 5000 & 276457 & 276457 & \textbf{0} \\
\hline[1pt]
\textbf{SA} & knapPI\_1\_10000\_1000\_1 & 10000 & 563647 & 563647 & \textbf{0} \\
GA & knapPI\_1\_10000\_1000\_1 & 10000 & 562556 & 563647 & 1091 \\
GSA & knapPI\_1\_10000\_1000\_1 & 10000 & 292255 & 563647 & 271392 \\
BB & knapPI\_1\_10000\_1000\_1 & 10000 & 130 & 563647 & 563517 \\
\textbf{CONGA} & knapPI\_1\_10000\_1000\_1 & 10000 & 563647 & 563647 & \textbf{0} \\
\hline[2pt]
\end{tblr}
\end{table}

\begin{table}[!htbp]
\centering
\caption{Results obtained on HD-WC datasets.}
\label{table:result-hd-wc}
\begin{tblr}{
  colspec = {llcccc}  
}
\hline[2pt]
  Algorithm &   Dataset &  Items & Best Value & Optimal Value & Shortfall   \\
\hline[2pt]
SA & knapPI\_2\_100\_1000\_1 & 100 & 1486 & 1514 & 28 \\
GA & knapPI\_2\_100\_1000\_1 & 100 & 1158 & 1514 & 356 \\
GSA & knapPI\_2\_100\_1000\_1 & 100 & 1041 & 1514 & 473 \\
BB & knapPI\_2\_100\_1000\_1 & 100 & 1440 & 1514 & 74 \\
\textbf{CONGA} & knapPI\_2\_100\_1000\_1 & 100 & 1514 & 1514 & \textbf{0} \\
\hline[1pt]
SA & knapPI\_2\_1000\_1000\_1 & 1000 & 6831 & 9052 & 2221 \\
GA & knapPI\_2\_1000\_1000\_1 & 1000 & 7912 & 9052 & 1140 \\
GSA & knapPI\_2\_1000\_1000\_1 & 1000 & 5675 & 9052 & 3377 \\
BB & knapPI\_2\_1000\_1000\_1 & 1000 & 9006 & 9052 & 46 \\
\textbf{CONGA} & knapPI\_2\_1000\_1000\_1 & 1000 & 9052 & 9052 & \textbf{0} \\
\hline[1pt]
SA & knapPI\_2\_10000\_1000\_1 & 10000 & 57852 & 90204 & 32352 \\
GA & knapPI\_2\_10000\_1000\_1 & 10000 & 79615 & 90204 & 10589 \\
GSA & knapPI\_2\_10000\_1000\_1 & 10000 & 54447 & 90204 & 35757 \\
BB & knapPI\_2\_10000\_1000\_1 & 10000 & 90073 & 90204 & 131 \\
\textbf{CONGA} & knapPI\_2\_10000\_1000\_1 & 10000 & 90204 & 90204 & \textbf{0} \\
\hline[1pt]
SA & knapPI\_2\_200\_1000\_1 & 200 & 1537 & 1634 & 97 \\
GA & knapPI\_2\_200\_1000\_1 & 200 & 1306 & 1634 & 328 \\
GSA & knapPI\_2\_200\_1000\_1 & 200 & 1073 & 1634 & 561 \\
BB & knapPI\_2\_200\_1000\_1 & 200 & 1603 & 1634 & 31 \\
\textbf{CONGA} & knapPI\_2\_200\_1000\_1 & 200 & 1634 & 1634 & \textbf{0} \\
\hline[1pt]
SA & knapPI\_2\_2000\_1000\_1 & 2000 & 12780 & 18051 & 5271 \\
GA & knapPI\_2\_2000\_1000\_1 & 2000 & 15887 & 18051 & 2164 \\
GSA & knapPI\_2\_2000\_1000\_1 & 2000 & 11064 & 18051 & 6987 \\
BB & knapPI\_2\_2000\_1000\_1 & 2000 & 17794 & 18051 & 257 \\
\textbf{CONGA} & knapPI\_2\_2000\_1000\_1 & 2000 & 18051 & 18051 & \textbf{0} \\
\hline[1pt]
SA & knapPI\_2\_500\_1000\_1 & 500 & 3744 & 4566 & 822 \\
GA & knapPI\_2\_500\_1000\_1 & 500 & 3701 & 4566 & 865 \\
GSA & knapPI\_2\_500\_1000\_1 & 500 & 2951 & 4566 & 1615 \\
BB & knapPI\_2\_500\_1000\_1 & 500 & 4484 & 4566 & 82 \\
\textbf{CONGA} & knapPI\_2\_500\_1000\_1 & 500 & 4566 & 4566 & \textbf{0} \\
\hline[1pt]
SA & knapPI\_2\_5000\_1000\_1 & 5000 & 29220 & 44356 & 15136 \\
GA & knapPI\_2\_5000\_1000\_1 & 5000 & 37746 & 44356 & 6610 \\
GSA & knapPI\_2\_5000\_1000\_1 & 5000 & 27387 & 44356 & 16969 \\
BB & knapPI\_2\_5000\_1000\_1 & 5000 & 44198 & 44356 & 158 \\
\textbf{CONGA} & knapPI\_2\_5000\_1000\_1 & 5000 & 44356 & 44356 & \textbf{0} \\
\hline[2pt]
\end{tblr}
\end{table}

\begin{table}[!htbp]
\centering
\caption{Results obtained on HD-SC datasets.}
\label{table:result-hd-sc}
\begin{tblr}{
  colspec = {llcccc}  
}
\hline[2pt]
  Algorithm &   Dataset &  Items & Best Value & Optimal Value & Shortfall \\
\hline[2pt]
SA & knapPI\_3\_100\_1000\_1 & 100 & 2296 & 2397 & 101 \\
GA & knapPI\_3\_100\_1000\_1 & 100 & 2091 & 2397 & 306 \\
GSA & knapPI\_3\_100\_1000\_1 & 100 & 1095 & 2397 & 1302 \\
BB & knapPI\_3\_100\_1000\_1 & 100 & 2268 & 2397 & 129 \\
\textbf{CONGA} & knapPI\_3\_100\_1000\_1 & 100 & 2397 & 2397 & \textbf{0} \\
\hline[1pt]
SA & knapPI\_3\_1000\_1000\_1 & 1000 & 11789 & 14390 & 2601 \\
GA & knapPI\_3\_1000\_1000\_1 & 1000 & 13090 & 14390 & 1300 \\
GSA & knapPI\_3\_1000\_1000\_1 & 1000 & 5589 & 14390 & 8801 \\
BB & knapPI\_3\_1000\_1000\_1 & 1000 & 14271 & 14390 & 119 \\
\textbf{CONGA} & knapPI\_3\_1000\_1000\_1 & 1000 & 14390 & 14390 & \textbf{0} \\
\hline[1pt]
SA & knapPI\_3\_10000\_1000\_1 & 10000 & 106114 & 146919 & 40805 \\
GA & knapPI\_3\_10000\_1000\_1 & 10000 & 124719 & 146919 & 22200 \\
GSA & knapPI\_3\_10000\_1000\_1 & 10000 & 54518 & 146919 & 92401 \\
BB & knapPI\_3\_10000\_1000\_1 & 10000 & 146787 & 146919 & 132 \\
\textbf{CONGA} & knapPI\_3\_10000\_1000\_1 & 10000 & 146919 & 146919 & \textbf{0} \\
\hline[1pt]
SA & knapPI\_3\_200\_1000\_1 & 200 & 2594 & 2697 & 103 \\
GA & knapPI\_3\_200\_1000\_1 & 200 & 25319 & 2697 & -22622 \\
GSA & knapPI\_3\_200\_1000\_1 & 200 & 1095 & 2697 & 1602 \\
BB & knapPI\_3\_200\_1000\_1 & 200 & 2542 & 2697 & 155 \\
\textbf{CONGA} & knapPI\_3\_200\_1000\_1 & 200 & 2697 & 2697 & \textbf{0} \\
\hline[1pt]
SA & knapPI\_3\_2000\_1000\_1 & 2000 & 22482 & 28919 & 6437 \\
GA & knapPI\_3\_2000\_1000\_1 & 2000 & 25319 & 28919 & 3600 \\
GSA & knapPI\_3\_2000\_1000\_1 & 2000 & 10818 & 28919 & 18101 \\
BB & knapPI\_3\_2000\_1000\_1 & 2000 & 28726 & 28919 & 193 \\
\textbf{CONGA} & knapPI\_3\_2000\_1000\_1 & 2000 & 28919 & 28919 & \textbf{0} \\
\hline[1pt]
SA & knapPI\_3\_500\_1000\_1 & 500 & 6103 & 7117 & 1014 \\
GA & knapPI\_3\_500\_1000\_1 & 500 & 6517 & 7117 & 600 \\
GSA & knapPI\_3\_500\_1000\_1 & 500 & 2916 & 7117 & 4201 \\
BB & knapPI\_3\_500\_1000\_1 & 500 & 6995 & 7117 & 122 \\
\textbf{CONGA} & knapPI\_3\_500\_1000\_1 & 500 & 7117 & 7117 & \textbf{0} \\
\hline[1pt]
SA & knapPI\_3\_5000\_1000\_1 & 5000 & 53672 & 72505 & 18833 \\
GA & knapPI\_3\_5000\_1000\_1 & 5000 & 61904 & 72505 & 10601 \\
GSA & knapPI\_3\_5000\_1000\_1 & 5000 & 27304 & 72505 & 45201 \\
BB & knapPI\_3\_5000\_1000\_1 & 5000 & 72345 & 72505 & 160 \\
\textbf{CONGA} & knapPI\_3\_5000\_1000\_1 & 500 & 72505 & 72505 & \textbf{0} \\
\hline[2pt]
\end{tblr}
\end{table}

\end{document}